\scrollmode
\documentclass[12pt]{amsart}
\usepackage{amssymb}
\usepackage{verbatim}
\usepackage{eucal}
\usepackage{graphicx}
\usepackage{psfrag}
\usepackage{mathrsfs}
\theoremstyle{plain}
\swapnumbers
\newtheorem{theorem}[subsection]{Theorem}
\newtheorem{proposition}[subsection]{Proposition}
\newtheorem{lemma}[subsection]{Lemma}

\theoremstyle{definition}
\newtheorem{definition}[subsection]{Definition}

\theoremstyle{remark}
\newcommand{\R}{{\mathbb R}}

\newcommand{\C}{{\mathbb C}}
\newcommand{\Z}{{\mathbb Z}}
\newcommand{\Q}{{\mathbb Q}}

\newcommand{\M}{\mathscr{M}}
\newcommand{\NNN}{\mathscr{N}}

\newcommand{\OOO}{\mathscr{O}}

\newcommand{\sign}{\operatorname{sgn}}

\DeclareMathOperator*{\vol}{vol}
\DeclareMathOperator*{\ord}{ord}
\DeclareMathOperator*{\diag}{diag}

\begin{document}

\title{Finiteness of minimal modular symbols for $SL_{n}$}

\author{Paul E. Gunnells}
\address{Department of Mathematics\\
Columbia University\\
New York, New York  10027}
\date{December, 1997. Revised May, 1999}
\email{gunnells@math.columbia.edu}
\subjclass{11F75, 11H06}
\keywords{Modular symbols, cohomology of arithmetic groups, 
geometry of numbers}

\thanks{The author would like to thank Avner Ash, Farshid Hajir, and
the referree for helpful comments.  The author was partially supported
by the NSF}

\begin{abstract}
Let $K/\Q $ be a number field with ring of integers $\OOO $,
and let $\Gamma \subset SL_{n} (\OOO )$ be a finite index subgroup.
Using a classical construction from the geometry of numbers and the
theory of modular symbols, we exhibit a finite spanning set of the
highest nonvanishing rational cohomology group of $\Gamma $.  
\end{abstract}

\maketitle

\section{Introduction}\label{intro}
Let $K/\Q $ be a number field with ring of integers $\OOO $.
Let $\Gamma \subset SL_{n} (\OOO )$ be a finite index subgroup, and let
$\nu $ be the virtual cohomological dimension of $\Gamma $.  That is, if
$\Gamma' \subset \Gamma $ is any finite index torsion-free subgroup,
then $H^{i} (\Gamma' ,M)=0$ for $i>\nu $ and any $\Z \Gamma $-module $M$.

Let $\M$ be the free abelian group generated by the symbols $[v_{1},\dots,
v_{n}]$, where the $v_{i}$ are nonzero points in $K^{n}$, modulo the
following relations:
\begin{enumerate}\label{mod.sym.def}
\item [(1)]  If $\tau$ is a
permutation on $n$ letters, then $[v_{1},\ldots,v_{n}] = \sign (\tau )
[\tau (v_{1}),\ldots,\tau (v_{n})]$, where $\sign (\tau)$ is  the sign
of $\tau $.
\item [(2)]  If $q\in K^{\times }$, then $[q v_{1},v_{2},\ldots,v_{n}] = [v_{1},\ldots,v_{n}]$.
\item [(3)]  If the $v_{i}$ are linearly dependent, then
$[v_{1},\ldots,v_{n}] = 0$.
\item [(4)]  If $v_{0},\ldots,v_{n}$ are nonzero points in $K^{n}$, then 
$\sum _{i}(-1)^{i}[v_{0},\ldots,\hat{v_{i}},\ldots,v_{n}]=0$.
\end{enumerate}
Elements of $\M$ are called \emph{minimal modular symbols}.  By a
theorem of Ash~\cite{ash.minmod} there is a surjective map of $\Z \Gamma
$-modules
\setcounter{equation}{4}
\begin{equation}\label{map.to.H}
\M \longrightarrow H^{\nu } (\Gamma ,\Q ).
\end{equation}
Hence $\M$ modulo $\Gamma $
is a spanning set for $H^{\nu } (\Gamma ,\Q )$.

However, this spanning set is not finite.  If $\OOO $ is a euclidean
domain, Ash and Rudolph~\cite{ash.rudolph} define a subset $\M_{u}\subset
\M$ with finite image under \eqref{map.to.H}, and give an efficient
algorithm to write any $[m]\in \M$ as a sum $[m]=\sum [m_{i}]$, where
each $[m_{i}]\in \M_{u}$.  Our goal in this note is more modest: in
Theorem~\ref{main} we exhibit a \emph{finite} spanning set for $H^{\nu}
(\Gamma ,\Q )$ for all $K$, not necessarily euclidean, but we provide
no practical reduction algorithm.  The proof relies on a classical
construction of Minkowski from the geometry of numbers.

\section{Statement of the result}\label{int.strut}
\subsection{}
Let $[K : \Q ]=d$.  For any ring $R$, let $M_{n} (R)$ be the $n\times
n$ matrices over $R$.

Given $m\in M_{n} (\OOO )$, let $[m]\in \M $ be the modular symbol
$[m_{1},\dots ,m_{n}]$, where the $m_{i}$ are the columns of $m$.  The
map $M_{n} (\OOO)\rightarrow \M $ is surjective, since using relations
(1) and (2) we have $[v_{1},\dots ,v_{n}] = [q_{1}v_{1},\dots
,q_{n}v_{n}]$ for any nonzero $v_{i}\in K^{n}$ and any $q_{i}\in
K^{\times }$.

Let $N_{K/\Q }\colon K\rightarrow \Q $ be the norm map.  Define a map
$\|\phantom{m}\|\colon M_{n} (K)\rightarrow \Q^{\geq 0}$ by 
\[
\|{m}\| = | N_{K/\Q } (\det m)|.
\]

\subsection{}
Let $C\geq 1$ be an integer, and define
\[
M (C) := \left\{m \in M_{n} (\OOO )\bigm| \|{m}\| \leq C \right\}.
\]
Let $\M (C)\subset \M $ be the set of modular symbols in the image of
$M (C)$ under the map $M_{n} (\OOO)\rightarrow \M $.

\begin{proposition}\label{finite}
Let $\Gamma \subset SL_{n} (\OOO )$ be of finite index.  Then for any
$C\geq 1$, the set $\Gamma
\backslash \M (C)$ is finite.
\end{proposition}

\begin{proof}
It suffices to verify the statement for $\Gamma =SL_{n} (\OOO )$.  To
simplify notation, we write $G$ for $GL_{n} (\OOO )$.  	If $m\in M_{n}
(\OOO )$, then we write $\Lambda (m)$ for the $\OOO$-lattice generated
by the columns of $m$.

First, we claim that $G\backslash M (C)$ is finite for any $C\geq 1$.
Indeed, the set 
\[
\bigl\{m\bigm | [\OOO ^{n}:\Lambda (m)]\leq C
\bigr\}
\]
is finite modulo $G$.  Now for any matrix $m$, the index
$[\OOO ^{n}:\Lambda (m)]$ is equal to $\NNN (\ord (T))$, where $\NNN$
denotes the ideal norm, and $\ord (T)$ is the \emph{order ideal} of
the torsion module $T:=\OOO ^{n}/\Lambda (m)$ (\cite{reiner}, \S4D).
Furthermore, $\ord (T)$ is a principal ideal generated by $\det (f)$,
where $f\colon \OOO ^{n}\rightarrow \OOO ^{n}$ is any $\OOO$-linear
endomorphism with image $\Lambda (m)$.  Clearly for $f$ we may take
multiplication by $m$, and thus
\begin{align*}
[\OOO ^{n}:\Lambda (m)] &= \NNN ((\det m))\\
	&= |N_{K/\Q } (\det m)|\\
	&=\|m\|.
\end{align*}
This implies $M (C)$ is finite modulo $G$.

We claim this implies $\M (C)$ is finite modulo $\Gamma $.
To see this we use the following easily verified fact.  Suppose a group
$A$ acts on a set $S$, and $B\vartriangleleft A$ is a normal
subgroup.  If $A/B$ is a finite group, and $A\backslash S$ is
finite, then $B\backslash S$ is finite.

To apply this, let $A=G$ and $B=Z\cdot \Gamma $, where $Z$ is the
center of $G$.  The group $G/ (Z\cdot \Gamma)$ is isomorphic to
$U/U^{n}$, where $U = \OOO ^{\times }$ and $U^{n}$ is the subgroup of
$n$th powers.  Hence $G/ (Z\cdot \Gamma)$ is finite by
Dirichlet's unit theorem.  Setting $S=M (C)$, we have that
$G\backslash M (C)$ finite implies that $(Z\cdot
\Gamma)\backslash M (C)$ is finite.  Since $M (C)$ maps surjectively
onto $\M (C)$, and $Z$ acts trivially on $\M$, the result follows.
\end{proof}

\subsection{}
Let $V=K\otimes \R \cong \R ^{r}\times \C ^{s}$, where $r+2s=d$.  We
define the \emph{Minkowski constant}\footnote{A classical result of
Minkowski bounds the discriminant $D_{K}$ of $K$: if $K\otimes \R
\cong \R ^{r}\times \C ^{s}$, then $D_{K}^{2}\geq M_{K}$.} $M_{K}$ by
\[
M_{K} = \left(\frac{\pi }{4} \right)^{s}\frac{d^{d}}{d!}.
\]

We now state our main result.
\begin{theorem}\label{main}
The image of $\M (C)$ under \eqref{map.to.H} spans $H^{\nu} (\Gamma ,\Q )$ for 
\[
C\geq \left\lfloor\left(\frac{\sqrt{|D_{K}|}}{M_{K}} \right)^{n}\right\rfloor.
\]
\end{theorem}

\section{Proof}\label{proof}
\subsection{}
Fix a modular symbol $[m]$, where $m\in M_{n} (\OOO )$.  We claim that
if $\| m \|> (\sqrt{|D_{K}|}/M_{K})^{n}$, then we can find a nonzero
$x\in \OOO ^{n}$ such that $x=\sum q_{i}v_{i}$ with $q_{i}\in K$ and
$|N_{K/\Q } (q_{i})|<1$.  We may then use (4) to construct a relation
$[m]=\sum (-1)^{i} [m_{i}]$, where
\[
[m_{i}]:= [v_{1},\dots
,v_{i-1},\hat v_{i},x,v_{i+1},\dots ,v_{n}].
\]
Clearly $\|m_{i}\| <
\|m\|$.  This claim implies the theorem, because by iterating this
process we can write any $[m]\in \M$ as a sum of symbols from $\M
(C)$.

To prove the claim we use the \emph{regular representation} of $\OOO
$.  Fix a $\Z $-basis $\omega _{1} = 1,\omega _{2},\dots ,\omega _{d}$
of $\OOO $.  Then this representation is the map $\OOO \rightarrow
M_{d} (\Z )$ defined by $\alpha \mapsto \ell _{\alpha }$, where $\ell
_{\alpha }$ is the matrix of the map $x\mapsto \alpha x$ in terms of
the $\Z $-basis.  This induces a ring homomorphism $\varphi \colon
M_{n} (\OOO )\rightarrow M_{nd} (\Z )$, in which matrix entries are
taken to $d\times d$ blocks.  Via $\varphi $, any column vector $v\in
\OOO ^{n}$ determines $d$ column vectors $\{v^{1},\dots, v^{d}
\}\subset \Z ^{nd}$.

Now apply $\varphi $ to $m$:
\[
(v_{1},\dots ,v_{n})\longmapsto (v_{1}^{1},\dots ,v_{1}^{d},\dots ,v_{n}^{1},\dots
,v_{n}^{d})\in M_{nd} (\Z ).
\]
The matrix $\varphi (m_{i})$ is obtained from $\varphi (m)$ by replacing the
columns $v_{i}^{1},\dots ,v_{i}^{d}$ with $x^{1},\dots ,x^{d}$.

For $1\leq i\leq n$, $1\leq j\leq d$, let $\lambda _{i}^{j}$ be
real variables, and consider the region $S\subset \R ^{nd}$ defined by 
\[
S := \Biggl\{\sum _{\substack{1\leq i\leq n\\
1\leq j\leq d}}\lambda _{i}^{j}v_{i}^{j}\Biggm | \Bigl|N_{K/\Q}
\Bigl(\sum_{j} \lambda _{i}^{j}\omega _{j}\Bigr)\Bigr| < 1,
\text{where $1\leq i\leq n$}\Biggr\}.
\]
Here we interpret $N_{K/\Q }(\sum_{j} \lambda _{i}^{j}\omega _{j})$ to mean
the polynomial in $\Z [\lambda _{i}^{j}]$ constructed using the
norm form.

\begin{lemma}\label{newlemma1}
Let $x\in \OOO $.  Suppose that $x=\sum _{i}q_{i}v_{i}$, where
$v_{i}\in \OOO $ and $q_{i}\in K$.  Write $q_{i}=\sum
_{j}q_{i}^{j}\omega _{j}$.  Then $x^{1} =
\sum_{i,j}q_{i}^{j}v_{i}^{j}$.
\end{lemma}

\begin{proof}
For any $x\in \OOO $, let $C_{k} (x)$ be the coefficient of $\omega
_{k}$ in the expansion of $x$ in terms of the fixed $\Z $-basis.  By
the definition of $\varphi $, the $k$th component of $x^{j}\in \Z ^{d}$ is
$C_{k} (\omega _{j}x)$.  In particular, since $\omega _{1}=1$, the
$k$th component of $x^{1}$ is $C_{k} (x) = C_{k} (\sum q_{i}v_{i})$.

Now let $y\in \Z ^{d}$ be the vector $\sum _{i,j}q^{j}_{i}v^{j}_{i}$.
We will show that the components of $y$ match those of $x^{1}$.
Indeed, the $k$th component of $y$ is   
$\sum _{i,j}q^{j}_{i}C_{k} (\omega
_{j}v_{i})$.  But then
\begin{align*}
\sum _{i,j}q^{j}_{i}C_{k}(\omega _{j}v_{i})&=\sum _{i,j}C_{k}(q^{j}_{i}\omega _{j}v_{i})\\
&=C_{k} (\sum _{i,j} q^{j}_{i}\omega _{j}v_{i})\\
&=C_{k} (\sum_{i} (\sum _{j}q^{j}_{i}\omega _{j})v_{i})\\
&=C_{k} (\sum _{i}q_{i}v_{i}).
\end{align*}
The final expression is the $k$th component of
$x^{1}$, so the result follows.
\end{proof}

\begin{lemma}\label{sameas}
There exists a nonzero $x\in \OOO ^{n}$ such that $x=\sum q_{i}v_{i}$
with $q_{i}\in K$ and $|N_{K/\Q } (q_{i})|<1$ if and only if the
region $S$ contains a nonzero rational integral point.
\end{lemma}

\begin{proof}
Let $x\in \OOO ^{n}$ satisfy the hypotheses.  Apply the regular
representation to $x=\sum q_{i}v_{i}$ and Lemma \ref{newlemma1} to
each row of $x$.  We find $x^{1} = \sum q^{j}_{i} v^{j}_{i}$, where
$q^{j}_{i}\in \Q $ and $q_{i}=\sum q_{i}^{j}\omega _{j}$.  The
condition $|N_{K/\Q } (q_{i})|<1$ is thus exactly $|N_{K/\Q }
(\sum_{j} q_{i}^{j}\omega _{j})| < 1$.  Hence $x^{1}$ is a nonzero
rational integral point in $S$.  The converse follows by reversing
this argument.
\end{proof}

\subsection{}
Now we will find a bounded symmetric convex body $P\subset S$ and show
that if $\| m \| >(\sqrt{|D_{K}|}/M_{K})^{n} $, then $\vol P >2^{nd}$.
Then by Minkowski's theorem (\cite{fr.tay}, IV.2.6)
$P$ and hence $S$ will contain a nonzero integral point.  By Lemma
\ref{sameas} this will imply Theorem \ref{main}.

To do this, apply $\varphi (m^{-1})$ to $S$.  This carries 
$\{v_{i}^{j}\}$ onto the standard basis of $\R ^{nd}$.  We can then write $\varphi
(m^{-1})(S)$ as the $n$-fold product $T^{n}$, where $T\subset \R
^{d}$ is the region
\[
T := \Biggl\{(y_{1},\dots ,y_{d})\Biggm | \Bigl|N_{K/\Q } \Bigl(\sum y_{i}\omega _{i}\Bigr)\Bigr|\leq
1 \Biggr\}.
\]
This region can be transformed further as follows.  The vector space $V$
contains $\OOO $, embedded by $\alpha \mapsto \alpha \otimes 1$.  Let
$\mu \colon \R ^{d}\rightarrow \R ^{d}$ be the linear map taking the
standard basis of $\Z ^{d}$ to $\{\omega _{1}\otimes 1,\dots ,\omega
_{d}\otimes 1 \}$.  Then
\[
\mu (T) = \left\{(x_{1},\dots ,x_{r},z_{1},\dots ,z_{s}) \in \R
^{r}\times \C ^{s}\bigm | |x_{1}\cdots x_{r}|z_{1}\bar z_{1}\cdots
z_{s}\bar z_{s} < 1 \right\}.
\] 

Now we construct the \emph{generalized octahedron} of Minkowski.
This will be a bounded, symmetric, convex
body in $\mu (T)$.  Take polar coordinates $(\rho
_{i},\theta _{i})$ for the $z_{i}$, and let $V^{+}\subset V$ be the
subset 
\[
V^{+} := \bigl\{(x_{1},\dots ,x_{r},\rho _{1},\theta _{1},\dots ,\rho
_{s},\theta _{s}) \bigm | \text{$x_{i}\geq 0$, $\rho _{i}\geq 0$, and $\theta _{i}=0$} \bigr\}.
\]

\begin{figure}[tbh]
\begin{center}
\psfrag{p}{$p$}
\psfrag{d}{$\Delta $}
\psfrag{mt}{$\mu (T)$}
\includegraphics{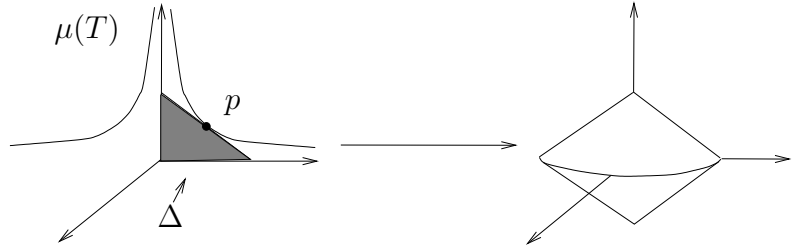}
\end{center}
\caption{\label{example}The generalized octahedron for $(r,s) = (1,1)$.}
\end{figure}

\begin{definition}\label{body}
Given a point $p\in \mu (T)\cap V^{+}$, let $Q (p)$ be the subset of $\mu (T)$ 
constructed as follows:
\begin{enumerate}
\item Construct the tangent hyperplane to $\mu (T)\cap V^{+}$ at $p$.
\item Use this hyperplane and the bounding hyperplanes of $V^{+}$ to
cut out a $( 2r)$-simplex $\Delta $ in $V^{+}$.
\item Apply the motions $x_{i}\mapsto -x_{i}$ and $\theta
_{i}\mapsto \theta _{i}+ \beta $, $0\leq \beta \leq 2\pi $ to $\Delta
$ to sweep out $Q (p)$ (cf. Figure~\ref{example}).
\end{enumerate}
\end{definition}

\begin{lemma}\label{facts}
$Q (p)$ is bounded, symmetric, and convex.  The volume of $Q (p)$ is independent
of $p$, and is 
\begin{equation}\label{volume}
2^{r+s}M_{K}.
\end{equation}
\end{lemma}

\begin{proof}
All statements are standard results from the geometry of
numbers~(\cite{fr.tay}, IV.2), except for the independence of $p$.  However,
this is easy to verify.
\end{proof}

\subsection{}
We now complete the proof of the theorem.  We choose $p\in \mu (T)\cap V^{+}$
and abbreviate $Q (p)$ to $Q$.  Define $P\subset S$ by 
\[
P := \varphi (m) \bigl( \mu ^{-1}Q \times \cdots \times \mu ^{-1}Q ).
\]
$P$ is symmetric, bounded, and convex, and 
\begin{equation}\label{volume2}
\vol P = |\det \varphi (m)|\left( \frac{\vol Q}{|\det \mu| }\right)^{n}.
\end{equation}
Now in \eqref{volume2} we apply Lemma \ref{facts}, and substitute
$|\det \mu | = 2^{-s}\sqrt{|D_{K}|}$ and 
$\| m \| = | \det \varphi (m) |$.  To ensure that $\vol P>2^{nd}$, we require
\[
\| m \|>\left(\frac{\sqrt{|D_{K}|}}{M_{K}} \right)^{n},
\]
as desired.

\bibliographystyle{amsplain}
\bibliography{bound}
\end{document}